\documentclass[review]{elsarticle}
\usepackage{lineno,hyperref}
\usepackage{amssymb}
\usepackage{mathrsfs}
\usepackage{color}
\usepackage{amsmath,amssymb,bm}
\usepackage{cases}
\modulolinenumbers[5]

\newtheorem{definition}{Definition}
\newtheorem{theorem}{Theorem}
\newtheorem{lemma}{Lemma}
\newtheorem{proposition}{Proposition}

\newtheorem{remark}{Remark}
\newtheorem{corollary}{Corollary}

\journal{arXiv}

\begin{document}

\begin{frontmatter}

\title{Existence of weak solutions for a class of non-divergent parabolic equations with variable exponent}

\author{Jingfeng Shao}
\ead{sjfmath@foxmail.com}
\author{Zhichang Guo\corref{mycorrespondingauthor}}
\cortext[mycorrespondingauthor]{Corresponding author}
\ead{mathgzc@hit.edu.cn}
\author{Zhongxiang Zhou}
\ead{zhouzx@hit.edu.cn}

\address{School of Mathematics, Harbin Institute of
Technology, Harbin, 150001, China}

\begin{abstract}
A doubly degenerate parabolic equation in non-divergent form with variable growth is investigated in this paper. In suitable spaces, we prove the existence of weak solutions of the equation for cases $1\leq m < 2$ and $m\geq 2$ in different ways. And we establish the non-expansion of support of the solution for the problem.
\end{abstract}


\end{frontmatter}


\section{Introduction}
Let $\Omega \subset \mathbb{R}^n$ be a bounded domain with Lipschitz
boundary $\partial \Omega$, and set $\Omega_T:=\Omega  \times (0,T)$, $\Gamma := \partial \Omega  \times (0,T)$. The goal of this article is to study the following diffusion problem:
\begin{numcases}{}
\frac{{\partial u}}{{\partial t}} = {u^m}{\rm div} {\left( {{{\left| {{D}u} \right|}^{{p}(x) - 2}}{D}u} \right)} \quad &\text{in } $\Omega_T$,\nonumber
 \\u(x,t) = 0 \quad &\text{on } $\Gamma$,\label{fu1}
 \\u(x,0) = u_0(x) \quad &\text{in } $\Omega$,\nonumber
\end{numcases}
 where $m\geq 1$, the variable exponent $p: \overline \Omega \rightarrow (1,\infty) $ is log-H\"older continuous functions, and $D=(D_1,D_2,\cdots,D_n)$, $D_i$ denotes the weak derivative with respect to $x_i$.

The problem \eqref{fu1} is a doubly degenerate parabolic equation with variable exponent in non-divergence form, which generalizes the evolutional $p(x)$-Laplace. Due to the degeneracy or singularity at $u = 0$ and $|Du| = 0$, the problem \eqref{fu1} does not have classical solution in general. In this paper, we only consider the non-negative weak solutions of the equation.

If $m<1$, we can transform the problem into a non-Newtonian polytropic filtration equation as follow
\begin{numcases}{}
\frac{{\partial v}}{{\partial t}} =  {{\rm div}\left( {{\left| {{D}\Psi (v)} \right|}^{{p}(x) - 2}{D}\Psi(v)} \right)}  \quad &\text{in } $\Omega_T$, \nonumber\\
  v(x,t) = 0 \quad &\text{on } $\Gamma$,\label{fte} \\
  v(x,0) = \Psi^{-1} ({u_0}) \quad &\text{in } $\Omega$, \nonumber
\end{numcases}
where
\begin{equation}
v = \Psi^{-1}(u) := \frac{{{u^{1 - m}}}}{{1 - m}},\;u = \Psi (v) := {\left( {(1 - m)v} \right)^{\frac{1}{{1 - m}}}}. \label{tran}
\end{equation}
The existence of stong solutions of this kind of equations have been investigated in \cite{MR2811763,MR3013415}. The blow-up and extinction of solutions have also been studied in some articles (see \cite{MR2438319,MR3130539}). In particular, if $m=0$, the problem becomes a parabolic $p(x)$-Laplace equation.

If $m\geq1$, the transform \eqref{tran} fails due to the equation has a lot of singularities at the boundary and inside ($v = +\infty$ when $u = 0$). But in this case, the equation \eqref{fu1} is equivalent to the following double degenerate parabolic equation in divergence type
\begin{equation}
\frac{{\partial u}}{{\partial t}} = {\rm div} {\left( {u^m}{{{\left| {{D}u} \right|}^{{p}(x) - 2}}{D}u} \right)}-mu^{m-1}|Du|^{p(x)}.
\end{equation}

For the case where $m\geq 1$ and $p(x)$ is a constant, there are some results on the equation \eqref{fu1} in a series of papers. In the case of  $p(x)\equiv 2$, Bertsch et al. investigate the non-uniqueness of solutions and some properties of viscosity solutions \cite{MR965742,MR1174811,MR1044287}, and Friedman \cite{MR853975} et al. study the blow-up of solutions. Such equations also appear in biological \cite{MR707172} or as models modelling the spread of an epidemic \cite{MR859613}.
In the case of $p(x)\equiv p \neq 2$, the problem has also been investigated during the past decades \cite{MR2902844,MR2113162,MR3635371}.

In our knowledge, when $m\geq 1$ and the exponent $p(x)$ is variable, there are few results. In recent years, we established the existence of weak solutions only for the case $1\leq m<2$ (see \cite{MR4107096}). That was because we can prove $Du_{n}$ (where $u_{n}$ represents the weak solution of the auxiliary equation) converge to $Du$ in $L^{p(x)}(\Omega_T)$ when $1\leq m<2$, but failed when $m\geq 2$. Therefore, the diffusion equations in non-divergence form still need to be studied. Nowadays, we have established the existence of weak solutions to the equation as we have found that $u_{n}^{\frac{m-1}{p(x)}}Du_{n}$ converge to $u^{\frac{m-1}{p(x)}}Du$ in $L^{p(x)}(\Omega_T)$ for the case of $m \geq 2$. It is worth mentioning that the uniqueness of the solution of the parabolic equation in non-divergence form \eqref{fu1} does not hold for $m\geq1$ in general (for example \cite{MR1174811,MR901093,MR2171907}).


The following existence theorem is the main results of this paper.
\begin{theorem}{\label{main theorem}}
Assume that $m\geq 1$, $0\leq u_0 \in L^\infty (\Omega)\cap W_0 ^{1,p(\cdot)}(\Omega)$, the problem \eqref{fu1} admits a weak solution.
\end{theorem}

This paper is organized as follows. In Section 2, we introduce some mathematical preliminaries. The Section 3 is devoted to the existence of weak solution of the problem. In Subsection 3.1, we consider an auxiliary problem and list some necessary results.
In Subsection 3.2, we prove the existence of weak solution. Finally in section 4, we investigate the non-expansion of support of the solution of the problem.

\section{Mathematical Preliminaries}
Set ${\Omega_\tau} = \Omega  \times (0,\tau]$ is a generic cylinder of an arbitrary finite height $\tau$. Throughout this paper, $(\cdot)_+$, $(\cdot)_-$ represent the cut-off functions, where $(s)_+:= \max \{s,0\}$, $(s)_-:= \min \{s,0\}$, $s\in \mathbb{R}$.
%
The following definitions of these function spaces
are based on \cite{MR1134951,MR2790542}.

We define the modular
\[{\varrho _{q( \cdot )}}(f): = \int_\Omega  | f(x){|^{q(x)}}dx.\]
Then the variable exponent Lebesgue space is defined as follows:
\[{L^{q( \cdot )}}(\Omega ): = \left\{ {u \text{ is measurable on } \Omega \text{ and satisfy } {{\varrho _{q( \cdot )}}(\lambda u) < \infty \text{ for some }\lambda  > 0} } \right\},\]
which is a Banach space equipped with the Luxemburg norm
$$
\|f\|_{q(\cdot),\Omega}:=\|f\|_{L^{q(\cdot)}(\Omega)}=\inf\left\{\alpha>0 \left| \varrho_{q(\cdot)}(f/\alpha)\leq
1\right. \right\}.
$$
If $q\in L^\infty(\Omega)$, define $q^- = \mathop {\mathrm{ess}\inf }\limits_{x \in \Omega } q(x)$, $q^+ = \mathop {\mathrm{ess}\sup }\limits_{x \in \Omega } q(x)$, and we denote by $q'(x)$ the conjugate exponent of $q(x)$ as follows:
$$q'(x) = \frac{{q(x)}}{{q(x) - 1}}.$$
In particular, for a bounded exponent, the following lemma holds (refer to \cite[Lemma 3.2.5]{MR2790542}).
\begin{lemma} \label{vels}
Let $q\in L^\infty (\Omega)$. For any $u\in L^{q(\cdot)}(\Omega)$ and ${\left\| u \right\|_{q( \cdot ),\Omega }} > 0$, we have
\begin{align}
\min \left\{ {{\varrho _{q( \cdot )}}{{(u)}^{\frac{1}{q^ - }}},{\varrho _{q( \cdot )}}{{(u)}^{\frac{1}{q^ + }}}} \right\} \leq {\varrho _{q( \cdot )}}(u) \leq \max \left\{ {{\varrho _{q( \cdot )}}{{(u)}^{\frac{1}{q^ - }}},{\varrho _{q( \cdot )}}{{(u)}^{\frac{1}{q^ + }}}} \right\}.
\end{align}
\end{lemma}

The Sobolev space $W^{1,q(\cdot)}(\Omega)$ is
defined by
$$
W^{1,q(\cdot)}(\Omega):=\left\{u\in L^{q(\cdot)}(\Omega)\left| |D u| \in
L^{q(\cdot)}(\Omega)\right.\right\},
$$
which is a Banach space equipped with the norm
$$
\|u\|_{W^{1,q(\cdot)}(\Omega)}:=\left\|Du\right\|_{q(\cdot),\Omega}+\|u\|_{q(\cdot),\Omega}.
$$
The space $W_0 ^{1,q(\cdot)}(\Omega)$ is the closure of $C^{\infty}_0(\Omega)$ (the set of smooth functions with compact support in $\Omega$) in the norm of $W^{1,q(\cdot)}(\Omega)$.

Assume that every $p$ is log-H\"older continuous and there exist constants $p^-$, $p^+$ such that
\begin{equation}\label{px}
 1<p^-\leq p(x) \leq
 p^+< \infty.
\end{equation}

We introduce the Banach space
\begin{align}
\mathcal {V}(\Omega ) = \left\{ {u(x) \left| u(x) \in {L^2}\left( \Omega  \right),{{\left| {{D}u(x)} \right|}^{p(x)}} \in {L^1}\left( \Omega  \right)\right.} \right\} \label{vspace}
\end{align}
with
\[{\left\| u \right\|_{\mathcal {V}(\Omega )}} = {\left\| u \right\|_{2,\Omega }} +  {{{\left\| {{D}u} \right\|}_{p( \cdot ),\Omega }}} ,\]
and the space $\mathcal {V}_0(\Omega )$ is the closure of $C^{\infty}_0(\Omega)$.

By $\mathcal{U}\left( {{\Omega _T}} \right)$ we denote the Banach space
\[\mathcal{U}\left( {{\Omega _T}} \right){\text{ = }}\left\{ {u:\left( {0,T} \right) \to \mathcal{V}\left( \Omega  \right)\left| {u \in {L^2}\left( {{\Omega _T}} \right),{{\left| {{D}u} \right|}^{p(x)}} \in {L^1}\left( {{\Omega _T}} \right)} \right.} \right\},\]
with
\[{\left\| u \right\|_{\mathcal{U}\left( {{\Omega _T}} \right)}} = {\left\| u \right\|_{2,{\Omega _T}}} +  {{{\left\| {{D}u} \right\|}_{( \cdot ),{\Omega _T}}}} .\]
We denote ${\mathcal{U}_0}\left( {{\Omega _T}} \right)$ as a subspace of $\mathcal{U}\left( {{\Omega _T}} \right)$ in which the elements have zero traces on $\Gamma$.
$\mathcal{U}'\left( {{\Omega _T}} \right)$ is the dual space (the space of bounded linear functionals) of $\mathcal{U}\left( {{\Omega _T}} \right)$ \cite{MR3328376}.
The norm in $\mathcal{U}'\left( {{\Omega _T}} \right)$ is defined by
\[{\left\| v \right\|_{\mathcal{U}'\left( {{\Omega _T}} \right)}} = \sup \left\{ {\left\langle {v,\varphi } \right\rangle \left| {\varphi  \in \mathcal{U}\left( {{\Omega _T}} \right),{{\left\| \varphi  \right\|}_{\mathcal{U}\left( {{\Omega _T}} \right)}}} \right.} \leq 1  \right \}.\]
\begin{definition}{\label{weaksense}}
A function $u(x,t)$ is called a weak solution of problem \eqref{fu1} provided that
\begin{itemize}
\item
$u \in \mathcal{U}\left( {{\Omega _T}} \right) \cap {L^\infty }\left( {{\Omega _T}} \right),\;\;\dfrac{{\partial u}}{{\partial t}} \in \mathcal{U}'\left( {{\Omega _T}} \right) \cap {L^2}\left( {{\Omega _T}} \right).$
\item
For every $\varphi \in C_0^1 (\Omega_T)$,
\begin{align}
\iint_{{\Omega _T}} {\frac{{\partial u}}{{\partial t}}}\varphi \mathrm{d}x\mathrm{d}t + \iint_{{\Omega _T}}  {{{\left| {{D}u} \right|}^{{p}(x) - 2}}{D}u \cdot {D}\left( {{u^m}\varphi } \right)} \mathrm{d}x\mathrm{d}t = 0.\label{weak form}
\end{align}
\item
The following equations hold in the sense of trace:
\begin{align}
u(x, t) = 0 \qquad &\text{on } \Gamma, \label{trace1}
\\
u(x,0) = u_0 (x) \qquad &\text{in } \Omega. \label{trace2}
\end{align}
\end{itemize}
\end{definition}

We recall also the following inequalities which are classical in the theory of $p$-Laplace equations. The proofs of the following lemmas are in the appendix.

\begin{lemma} \label{pmin}
For all $\xi$, $\eta \in \mathbb{R}^n$, the following inequalities hold:
\begin{itemize}
\item [\rm{(i)}]
If  $\;2\leq p< \infty$, $\left( {{{\left| \xi  \right|}^{p - 2}}\xi  - {{\left| \eta  \right|}^{p - 2}}\eta } \right) \cdot \left( {\xi  - \eta } \right) \geq  \frac{1}{{{2^{p - 1}}}}{\left| {\xi  - \eta } \right|^p}$;
\item [\rm{(ii)}]
If  $\;1\leq p<2$, $\left( {{{\left| \xi  \right|}^{p - 2}}\xi  - {{\left| \eta  \right|}^{p - 2}}\eta } \right) \cdot \left( {\xi  - \eta } \right) \geq (p - 1){\left( {{{\left| \xi  \right|}^p} + {{\left| \eta  \right|}^p}} \right)^{\frac{{p - 2}}{p}}}{\left| {\xi  - \eta } \right|^2}$.
\end{itemize}
\end{lemma}

A generalized H\"{o}lder's inequality is stated in the following lemma.
\begin{lemma}\label{holder}
Assume that $q(x):\Omega\to[1,+\infty)$ is a measurable function. For every $f\in L^{q(\cdot)}(\Omega)$ and $g\in L^{q'(\cdot)}(\Omega)$ the following inequality holds:
\begin{align}
\int_\Omega  {\left| {fg} \right|} dx \leq 2{\left\| f \right\|_{q( \cdot ),\Omega }}{\left\| g \right\|_{q'( \cdot ),\Omega }}.
\end{align}
\end{lemma}

\begin{lemma} \label{lest}
Let $p(x)$ is a measurable function such that $1<p^-\leq p(x) \leq p^+\leq 2$. Suppose that $Du,Dv\in L^{p(\cdot)}(\Omega)$ and ${\left\| Du \right\|_{p( \cdot ),\Omega }}+{\left\| Dv \right\|_{p( \cdot ),\Omega }}\neq 0$. Then
\begin{align}
&\int_\Omega  {\left( {{{\left| Du \right|}^{p(x) - 2}}Du - {{\left| Dv \right|}^{p(x) - 2}}Dv} \right) \cdot \left( {Du - Dv} \right)\mathrm{d}x} \nonumber\\
& \geqslant ({p^ - } - 1) {\left( {\frac{{\int_\Omega  {{{\left| {Du - Dv} \right|}^{p(x)}}\mathrm{d}x} }}{{2{{\left\| {{{(|Du{|^{p( \cdot )}} + |Dv{|^{p( \cdot )}})}^{\frac{{2 - p( \cdot )}}{2}}}} \right\|}_{\frac{2}{{2 - p( \cdot )}},\Omega }}}}} \right)^\lambda},
\end{align}
where $\lambda \in \{\frac{2}{p^-},\frac{2}{p^+}\}$.
\end{lemma}

\begin{remark} \label{rest}
If $\;2\leq p^- \leq p^+< \infty$, using $\rm{(i)}$ of {\rm Lemma \ref{pmin}}, one has
\[\int_\Omega  {\left( {{{\left| Du \right|}^{p(x) - 2}}Du - {{\left| Dv \right|}^{p(x) - 2}}Dv} \right) \cdot \left( {Du - Dv} \right)\mathrm{d}x}  \geqslant \frac{1}{{{2^{{p^ + } - 1}}}}\int_\Omega  {|Du - Dv{|^{p(x)}} \mathrm{d}x} .\]
\end{remark}

\section{The Existence of Weak Solution}
\subsection{The Regularized Problem and Auxiliary Results}
In this subsection, we employ the regularization method and obtain some auxiliary results to prove the existence of weak solution to the problem \eqref{fu1}.

Now we consider the following regularized problem:
\begin{numcases}{}
\dfrac{{\partial u}}{{\partial t}} = {u^m}{\rm div} {\left( {{{\left| {{D}u} \right|}^{{p}(x) - 2}}{D}u} \right)} \quad &\text{in } $\Omega_T$, \nonumber
\\
u(x, t) =\varepsilon \quad &\text{on } $\Gamma$,\label{reg1}
\\
u(x,0) = u_0 (x)+\varepsilon\quad &\text{in } $\Omega$, \nonumber
\end{numcases}
where $u_0 \in L^{\infty}(\Omega)$ and $u_0 \geq0$.

\begin{definition}
A function $u(x,t)$ is called a weak solution of regularized problem \eqref{reg1} provided that
\begin{itemize}
\item
$u \in \mathcal{U}\left( {{\Omega _T}} \right) \cap {L^\infty }\left( {{\Omega _T}} \right),\;\;\dfrac{{\partial u}}{{\partial t}} \in \mathcal{U}'\left( {{\Omega _T}} \right).$
\item
For every $\varphi \in C_0^1 (\Omega_T)$,
\begin{align}
\iint_{{\Omega _T}} {\frac{{\partial u}}{{\partial t}}}\varphi \mathrm{d}x\mathrm{d}t + \iint_{{\Omega _T}} {{{\left| {{D}u} \right|}^{{p}(x) - 2}}{D}u \cdot {D}\left( {{u^m}\varphi } \right)} \mathrm{d}x\mathrm{d}t = 0.\label{reg weak form}
\end{align}
\item
The following equations hold in the sense of trace:
\begin{align}
u(x, t) = \varepsilon \qquad &\text{on } \Gamma,\label{trace3}
\\
u(x,0) = u_0 (x)+\varepsilon  \qquad &\text{in } \Omega.\label{trace4}
\end{align}
\end{itemize}
\end{definition}

The regularized problem \eqref{reg1} is still in non-divergent. Through nonlinear transformation, we can transform it into a divergent diffusion equation and obtain the following three propositions (The proof is in the appendix).
\begin{proposition}{\label{reex}}
Assume that $m>0$ and $p(x)$ is log-H\"older continuous which satisfies \eqref{px}, $u_0 \in L^{\infty}(\Omega)$ and $u_0 \geq0$, then the problem \eqref{reg1} admits a weak solution.
\end{proposition}

Denote the solution of the regularized problem \eqref{reg1} as $u_\varepsilon$ with the parameter $\varepsilon$.
\begin{proposition}{\label{bj}}
Let the conditions in Proposition \ref{reex} be fulfilled, and assume $\;0<\varepsilon_1\leq\varepsilon_2$, then we have $u_{\varepsilon_1}\leq u_{\varepsilon_2}$.
\end{proposition}

\begin{proposition}{\label{Estimates}}
Let the conditions in Proposition \ref{reex} be fulfilled, and assume $0\leq u_0 \in L^\infty (\Omega)\cap W_0 ^{1,q(\cdot)}(\Omega)$. Then
\begin{align}
&\iint_{{\Omega _T}} {u_\varepsilon ^{ - m}}{\left( {\frac{{\partial {u_\varepsilon }}}{{\partial t}}} \right)^2}\mathrm{d}x\mathrm{d}t \leq C, \\
&\mathop{\rm ess\;sup}\limits_{t\in (0,T)} {\int_\Omega  {\frac{1}{{{p}(x)}}{{\left| {{D}{u_\varepsilon (x,t) }} \right|}^{p(x)}}} \mathrm{d}x}\leq C,\\
&0< \varepsilon \leq \mathop{\rm ess\;inf}\limits_{(x,t)\in \Omega_T} u_\varepsilon \leq \mathop{\rm ess\;sup}\limits_{(x,t)\in \Omega_T} u_\varepsilon \leq \mathop{\rm ess\;sup}\limits_{x\in \Omega} u_0+\varepsilon \leq C,
\end{align}
where $C$ is a constant independent of $\varepsilon$ and $T$.
\end{proposition}


\subsection{The Existence of the Weak Solution for the Problem \eqref{fu1}}
In this subsection, we devote to prove the existence of weak solutions to the equation. For $1\leq m<2$, we can get that there is a subsequence $u_n$ of $\left\{ {{u_\varepsilon }} \right\}$ converges to $u$ in $L^{p(x)}(\Omega_T)$. Although this result cannot be obtained for the case $m\geq 2$, as an alternative, we obtain that $u_n ^{m-1}|Du_n|^{p(x)}$ converges to $u^{m-1}|Du|^{p(x)}$ in $L^{1}(\Omega_T)$. That implies there is a subsequence which satisfies $u_{n}^{\frac{m-1}{p(x)}}Du_{n}$ converges to $u^{\frac{m-1}{p(x)}}Du$ in $L^{p(x)}(\Omega_T)$.

{\bf Proof of Theorem \ref{main theorem}}:

Based on the estimates of Proposition \ref{Estimates}, we can extract from $\left\{ {{u_\varepsilon }} \right\}$, a subsequence (labeled $\left\{ {{u_n}} \right\}$) such that
\begin{align}
&{u_n} \to u\;\;\text{in}\;\;{L^r}({\Omega_T}),\;\;\text{$r>0$ and a.e. in $\Omega_T$}, \label{cv1}
\\
&\dfrac{{\partial {u_n}}}{{\partial t}}\rightharpoonup\dfrac{{\partial u}}{{\partial t}}\;\;\text{in}\;\;{L^{2}}(\Omega_T),\label{cv2}\\
&D{u_n } \rightharpoonup D u\;\;\text{in}\;\;L^{p(\cdot)}(\Omega_T). \label{cv3}
\end{align}
Notice that
\begin{align}
\iint_{{\Omega _T}} {\frac{{\partial {u_n }}}{{\partial t}}\phi }\mathrm{d}x\mathrm{d}t + \iint_{{\Omega _T}} { {{{\left| {{D}{u_n }} \right|}^{{p}(x) - 2}}{D}{u_n } \cdot {D}\left( {u_n ^m\phi } \right)} }\mathrm{d}x\mathrm{d}t = 0,\label{fu4}
\end{align}
$\text{for any }\phi (x,t) \in C_0^1\left( {{\Omega_T}} \right).$ Since $C_0^1\left( {{\Omega_T}} \right)$ is dense in $L^{p^-}(0,T;\mathcal{V}_0(\Omega))$,
we choose
\[\phi (x,t) = u_n ^{ - k}\left( {{u_n } - \varepsilon_n  - u} \right),\quad 1\leq k.\]
Then we have
\begin{align}
&\iint_{{\Omega_T}} {\frac{{\partial {u_n }}}{{\partial t}}u_n ^{ - k}\left( {{u_n } - \varepsilon_n  - u} \right)}\mathrm{d}x\mathrm{d}t \nonumber
\\+& \iint_{{\Omega _T}} {{\left| {{D}{u_n}} \right|}^{{p}(x) - 2}{D}{u_n } \cdot {D}\left(u_n^{m-k} ({{u_n } - \varepsilon_n  - u}) \right) }\mathrm{d}x\mathrm{d}t = 0. \label{34k}
\end{align}

\textbf{For the case $1\leq m<2$}. We choose $k=m$,
note that
\begin{align*}
  &\iint_{{\Omega_T}} {\frac{{\partial {u_n }}}{{\partial t}}u_n ^{ - m}\left( {{u_n } - \varepsilon_n  - u} \right)}\mathrm{d}x\mathrm{d}t  \\
 = &\iint_{{\Omega_T}} {\varepsilon_n \frac{{\partial {u_n }}}{{\partial t}}u_n ^{ - m} }\mathrm{d}x\mathrm{d}t + \iint_{{\Omega_T}} {\frac{{\partial {u_n }}}{{\partial t}}u_n ^{ - m}\left( {u_n - u} \right)}\mathrm{d}x\mathrm{d}t \\
 =:&\rm{I}+\rm{II}.
\end{align*}
Proposition \ref{bj} implies that $u_\varepsilon$ converges to $u$ monotonically.
We have
\begin{align*}
  \left| \rm{I} \right|& \leqslant {\left( {\iint_{{\Omega_T}} {u_n ^{ - m}}{{\left( {\frac{{\partial {u_n }}}{{\partial t}}} \right)}^2}\mathrm{d}x\mathrm{d}t} \right)^{\frac{1}{2}}}{\left( {\iint_{{\Omega_T}} {u_n ^{m-2k}}{\varepsilon_n ^2}\mathrm{d}x\mathrm{d}t} \right)^{\frac{1}{2}}} \\
  & \leqslant C{\varepsilon_n ^{\frac{{2 - m}}{2}}}{\left( {\iint_{{\Omega_T}} {\frac{{{\varepsilon_n ^m}}}{{u_n ^m}}}\mathrm{d}x\mathrm{d}t} \right)^{\frac{1}{2}}},
\end{align*}
\begin{align*}
\left| {\rm{II}} \right| \leqslant {\left( {\iint_{{\Omega_T}} {u_n ^{ - m}}{{\left( {\frac{{\partial {u_n }}}{{\partial t}}} \right)}^2}\mathrm{d}x\mathrm{d}t} \right)^{\frac{1}{2}}}{\left( {\iint_{{\Omega_T}} {\left( {u_n ^{1 - \frac{m}{2}} - \frac{u}{{u_n ^{\frac{m}{2}}}}} \right)}\mathrm{d}x\mathrm{d}t} \right)^{\frac{1}{2}}}.
\end{align*}
Using Lebesgue's dominated convergence theorem, we get  \[\iint_{{\Omega_T}} {\frac{{\partial {u_n }}}{{\partial t}}u_n ^{ - m}\left( {{u_n } - \varepsilon_n  - u} \right)}\mathrm{d}x\mathrm{d}t \to 0,\;\text{as $\varepsilon_n \rightarrow 0$}.\]
Therefore, by \eqref{34k}, we have
\[\iint_{{\Omega _T}} {{{{\left| {{D}{u_n }} \right|}^{{p}(x) - 2}}{D}{u_n } \cdot {D}\left( {{u_n } - u} \right)} }\mathrm{d}x\mathrm{d}t \to 0,\;\text{when}\; \varepsilon_n \rightarrow 0.\]
With \eqref{cv3},
we deduce
\begin{align}
 {\iint_{{\Omega _T}} {\left( {{{\left| {{D}{u_n }} \right|}^{{p}(x) - 2}}{D}{u_n } - {{\left| {{D}u} \right|}^{{p}(x) - 2}}{D}u} \right) \cdot {D}\left( {{u_n} - u} \right)}\mathrm{d}x\mathrm{d}t}  \to 0.\label{cv4}
\end{align}
Note
\[\Omega _{T}^1 := \left\{ {(x,t) \in {\Omega _T}\left| {\;1 < {p}(x) < 2} \right.} \right\},\]
\[\Omega _{T}^2 := {\Omega _T}\backslash \Omega _{T}^1.\]
For any fixed $x$, we have
\[\left( {{{\left| {{D}{u_n }} \right|}^{{p}(x) - 2}}{D}{u_n } - {{\left| {{D}u} \right|}^{{p}(x) - 2}}{D}u} \right) \cdot \left( {{D}{u_n } - {D}u} \right) \geqslant 0.\]
Furthermore,
\begin{align}
  \iint_{\Omega _{T}^1} {\left( {{{\left| {{D}{u_n }} \right|}^{{p}(x) - 2}}{D}{u_n } - {{\left| {{D}u} \right|}^{{p}(x) - 2}}{D}u} \right) \cdot {D}\left( {{u_n } - u} \right)}\mathrm{d}x\mathrm{d}t \to 0, \label{cv5}\\
  \iint_{\Omega _{T}^2} {\left( {{{\left| {{D}{u_n }} \right|}^{{p}(x) - 2}}{D}{u_n } - {{\left| {{D}u} \right|}^{{p}(x) - 2}}{D}u} \right) \cdot {D}\left( {u_n - u} \right)}\mathrm{d}x\mathrm{d}t \to 0. \label{cv6}
\end{align}
According to Lemma \ref{lest} and \eqref{cv5}, it follows that
\[\iint_{\Omega _{T}^1} {{{\left| {{D}\left( {u_n - u} \right)} \right|}^{p(x)}}}\mathrm{d}x\mathrm{d}t \to 0.\]
Similarly, by Remark \ref{rest} and \eqref{cv6}, we have
\[\iint_{\Omega _{T}^2} {{{\left| {{D}\left( {u_n - u} \right)} \right|}^{p(x)}}}\mathrm{d}x\mathrm{d}t \to 0.\]
Therefore, we obtain
\[\iint_{\Omega _{T}} {{{\left| {{D}\left( {u_n - u} \right)} \right|}^{p(x)}}}\mathrm{d}x\mathrm{d}t \to 0.\]
which implies that
\begin{align}
D {u_n } \to D u \quad \text{in}\; L^{p(\cdot)}(\Omega_T). \label{cv7}
\end{align}

From \eqref{fu4}, we observe that
\begin{align*}
 &\iint_{{\Omega _T}} {\frac{{\partial {u_n }}}{{\partial t}}\phi }\mathrm{d}x\mathrm{d}t + {\iint_{{\Omega _T}} {m u_n ^{m - 1}{{\left| {{D}{u_n }} \right|}^{{p}(x)}}\phi }\mathrm{d}x\mathrm{d}t}
 \\+& {\iint_{{\Omega _T}} {u_n ^m{{\left| {{D}{u_n }} \right|}^{p(x) - 2}}{D}{u_n } \cdot {D}\phi }\mathrm{d}x\mathrm{d}t}  = 0.
\end{align*}
Combining with \eqref{cv1}, \eqref{cv2}, \eqref{cv7} and using Lebesgue's dominated convergence theorem, then we have
\begin{align*}
 &\iint_{{\Omega _T}} {\frac{{\partial {u }}}{{\partial t}}\phi }\mathrm{d}x\mathrm{d}t + {\iint_{{\Omega _T}} {mu ^{m - 1}{{\left| {{D}{u }} \right|}^{{p}(x)}}\phi }\mathrm{d}x\mathrm{d}t}
 \\+& {\iint_{{\Omega _T}} {{u }^m{{\left| {{D}{u }} \right|}^{p(x) - 2}}{D}{u } \cdot {D}\phi }\mathrm{d}x\mathrm{d}t}  = 0.
\end{align*}
Considering the limiting process, we have \eqref{trace1} and \eqref{trace2} in the sense of trace.

\textbf{For the case $m\geq2$.} Choose $k=1$ in \eqref{34k}, similar to the case $1\leq m< 2$, we get
\begin{align*}
\iint_{{\Omega_T}} {\frac{{\partial {u_n }}}{{\partial t}}u_n ^{ - 1}\left( {{u_n } - \varepsilon_n  - u} \right)}\mathrm{d}x\mathrm{d}t \to 0, \quad\text{when}\; \varepsilon_n \to 0.
\end{align*}
Thus,
\begin{align*}
\iint_{\Omega_T}{|Du_n|^{p(x)-2}Du_n}\cdot D\left(u_n ^{m-1}( u_n -u)\right)\mathrm{d}x\mathrm{d}t\to 0,\quad \text{when}\;\varepsilon_n \to 0.
\end{align*}
Then
\begin{align}
&\iint_{\Omega_T} {(m-1)u_n^{m-2}(u_n-\varepsilon_n-u) } |Du_n |^{p(x)}\mathrm{d}x\mathrm{d}t+\nonumber\\
&\iint_{\Omega_T}{u_n^{m-1}|Du_n|^{p(x)-2}Du_n }\cdot D( u_n-u)\mathrm{d}x\mathrm{d}t\to 0,\quad \text{when}\;\varepsilon_n \to 0.\label{341}
\end{align}
Since $u_n \in L^\infty(\Omega)$, $|Du_n|\in L^{p(\cdot)}(\Omega)$ and \eqref{cv1}, we have
\begin{align*}
& u_n ^{m-1}|Du|^{p(x)-2}Du \to u^{m-1}|Du|^{p(x)-2}Du \quad\text{in}\; L^{p'(x)}(\Omega_T),\\
& Du_n \rightharpoonup Du \quad\text{in}\; L^{p(x)}(\Omega_T).
\end{align*}
Thus,
\begin{align}
\iint_{\Omega_T}{u_n^{m-1}|Du|^{p(x)-2}Du}\cdot D( u_n-u)\mathrm{d}x\mathrm{d}t\to 0,\quad \text{when}\;\varepsilon_n \to 0. \label{342}
\end{align}
Combining with \eqref{341} and \eqref{342}, we get
\begin{align*}
&\iint_{\Omega_T} {u_n^{m-1}}\left({|Du_n|^{p(x)-2}Du_n}-|Du|^{p(x)-2}Du\right)\cdot D( u_n-u)\mathrm{d}x\mathrm{d}t+\nonumber\\
&\iint_{\Omega_T} {(m-1)u_n ^{m-2}(u_n-\varepsilon_n-u) } |Du_n|^{p(x)}\mathrm{d}x\mathrm{d}t \to 0,\quad \text{when}\;\varepsilon_n \to 0.
\end{align*}
Considering that
\begin{align*}
&\iint_{\Omega_T} {u_n^{m-1}}\left({|Du_n|^{p(x)-2}Du_n}-|Du|^{p(x)-2}Du\right)\cdot D( u_n-u)\mathrm{d}x\mathrm{d}t\ge 0,\\
&\iint_{\Omega_T} {(m-1)u_n ^{m-2}(u_n-u) } |Du_n|^{p(x)}\mathrm{d}x\mathrm{d}t \ge  0,\\
&\varepsilon_n \iint_{\Omega_T} {(m-1)u_n ^{m-2}} |Du_n|^{p(x)}\mathrm{d}x\mathrm{d}t \to 0,\quad \text{when}\;\varepsilon_n \to 0.
\end{align*}
Therefore,
\begin{align}
&\lim\limits_{\varepsilon_n\to 0_+} \iint_{\Omega_T} {u_n^{m-1}}\left({|Du_n|^{p(x)-2}Du_n}-|Du|^{p(x)-2}Du\right)\cdot D( u_n-u)\mathrm{d}x\mathrm{d}t =0,\label{343}\\
&\lim\limits_{\varepsilon_n \to 0_+} \iint_{\Omega_T} {(m-1)u_n ^{m-2}(u_n-u) } |Du_n|^{p(x)}\mathrm{d}x\mathrm{d}t = 0 . \label{344}
\end{align}
According to Proposition \ref{bj} and \eqref{343}, we have
\begin{align*}
\lim\limits_{\varepsilon_n\to 0_+} \iint_{\Omega_T} {u^{m-1}}\left({|Du_n|^{p(x)-2}Du_n}-|Du|^{p(x)-2}Du\right)\cdot D( u_n-u)\mathrm{d}x\mathrm{d}t =0.
\end{align*}
Thus, by Lemma \ref{lest} and Remark \ref{rest}, it follows that
\begin{align*}
u^{\frac{m-1}{p(x)}}Du_n \to u^{\frac{m-1}{p(x)}}Du \quad \text{in}\; L^{p(x)}(\Omega_T),
\end{align*}
which implies that (for a subsequence of $\{u_n\}$ if necessary, still labeled $\{u_n\}$)
\begin{align}
& u^{m-1}|Du_n|^{p(x)}\to u^{m-1}|Du|^{p(x)} \quad\text{in}\; L^{1}(\Omega_T),\label{345}\\
&u^{\frac{m-1}{p'(x)}}|Du_n |^{p(x)-2}Du_n \to u^{\frac{m-1}{p'(x)}}|Du|^{p(x)-2} Du \quad\text{in}\; L^{p'(\cdot)}(\Omega_T). \label{346}
\end{align}

We claim that
\begin{align}
& u_n ^{m-1}|Du_n|^{p(x)}\to u^{m-1}|Du|^{p(x)} \quad\text{in}\; L^{1}(\Omega_T),\label{347}\\
&u_n ^{m}|Du_n|^{p(x)-2}Du_n \to u^{m}|Du|^{p(x)-2} Du \quad\text{in}\; L^{1}(\Omega_T). \label{348}
\end{align}

Since
\begin{align*}
&\iint_{\Omega_T}\left| u_n ^{m-1}|Du_n|^{p(x)}-u^{m-1}|Du|^{p(x)}\right|\mathrm{d}x\mathrm{d}t\\
=&\iint_{\Omega_T} (u_n ^{m-1}-u^{m-1})|Du_n|^{p(x)}\mathrm{d}x\mathrm{d}t+\iint_{\Omega_T}\left| u ^{m-1}|Du_n|^{p(x)}-u^{m-1}|Du|^{p(x)}\right|\mathrm{d}x\mathrm{d}t\\
=:&A_1+A_2.
\end{align*}
By \eqref{344} and the Lagrange's mean value theorem, there exists $u\leq\xi\leq u_n$ satisfying
\begin{align*}
A_1=&\iint_{\Omega_T}{(m-1)(u_n-u) }\xi^{m-2} |Du_n|^{p(x)}\mathrm{d}x\mathrm{d}t\\
\leq & \iint_{\Omega_T}{(m-1)(u_n-u) }u_n^{m-2} |Du_n|^{p(x)}\mathrm{d}x\mathrm{d}t\to 0.
\end{align*}
By \eqref{345}, one has $A_2\to 0$.
Thus, we conclude \eqref{347}.

On the other hand, since
\begin{align*}
&\iint_{\Omega_T}\left| u_n ^{m}|Du_n|^{p(x)-2}Du_n -u^{m}|Du|^{p(x)-2} Du\right|\mathrm{d}x\mathrm{d}t\\
\leq &\iint_{\Omega_T}(u_n ^{m}-u^{m})|Du_n|^{p(x)-1}\mathrm{d}x\mathrm{d}t+\iint_{\Omega_T}{u^{m}}\left| |Du_n|^{p(x)-2}Du_n -|Du|^{p(x)-2} Du \right|\mathrm{d}x\mathrm{d}t\\
=:&B_1+B_2.
\end{align*}
By the  generalized H\"{o}lder's inequality, it follows that
\begin{align*}
B_1\leq & 2\left\||Du_n|^{p(x)-1}\right\|_{L^{p'(x)}(\Omega)}\cdot\|u_n ^m -u^m\|_{L^{p(x)}(\Omega)}\\
\leq & C \|u_n  -u\|_{L^{p(x)}(\Omega)} \to 0.
\end{align*}
By \eqref{345} and the continue embedding of $L^{p'(\cdot)}(\Omega)\subset L^{1}(\Omega)$, one has
\begin{align*}
& u_n^{\frac{m-1}{p'(x)}}|Du_n|^{p(x)-2}Du_n \to u^{\frac{m-1}{p'(x)}}|Du|^{p(x)-2}Du \quad \text{in}\; L^{1}(\Omega_T),
\end{align*}
then
\begin{align*}
B_2 = &\iint_{\Omega_T}{u^{1+\frac{m-1}{p(x)}}u^{\frac{m-1}{p'(x)}}}\left| |Du_n|^{p(x)-2}Du_n -|Du|^{p(x)-2} Du \right|\mathrm{d}x\mathrm{d}t\\
\leq &C \iint_{\Omega_T} \left| u^{\frac{m-1}{p'(x)}}|Du_n|^{p(x)-2}Du_n -u^{\frac{m-1}{p'(x)}}|Du|^{p(x)-2} Du \right|\mathrm{d}x\mathrm{d}t\to 0.
\end{align*}
Therefore, we conclude \eqref{348}.

Letting $\varepsilon_n \to 0$, combining \eqref{cv2}, \eqref{347} and \eqref{348}, we arrive at
\begin{align*}
 &\iint_{{\Omega _T}} {\frac{{\partial {u }}}{{\partial t}}\phi }dxdt + {\iint_{{\Omega _T}} {mu ^{m - 1}{{\left| {{D}{u }} \right|}^{{p}(x)}}\phi }\mathrm{d}x\mathrm{d}t}
 \\+& {\iint_{{\Omega _T}} {{u }^m{{\left| {{D}{u }} \right|}^{p(x) - 2}}{D}{u } \cdot {D}\phi }\mathrm{d}x\mathrm{d}t}  = 0.
\end{align*}
Considering the limiting process, we have \eqref{trace1} and \eqref{trace2} in the sense of trace.$\hfill\blacksquare$


\section{The Non-expansion of Support of the Solution}
For the case $m\geq1$, the solution has a localization property of non-expansion of the support. For a function $f:\Omega  \to {\mathbb{R}_ + } \cup \{ 0\}$, we denote by
$F$ the set $\{x\in \Omega| f(x)>0\}$ and define
\[\operatorname{supp} \;f: = \overline {\left\{ {x \in F\left| {\mathop {\lim }\limits_{r \to 0} \frac{{\mu \left( {F \cap {B_r}(x)} \right)}}{{\mu \left( {{B_r}(x)} \right)}} > 0} \right.} \right\}},\]
where ${B_r}(x)$ denotes $\left\{ {z\left|\;{\left| {z - x} \right| \leqslant r} \right.} \right\}$.
\begin{theorem} \label{nochange}
Let u be a weak solution of \eqref{fu1} with \rm{supp} $u_0 \subsetneqq \Omega$, $u_0\geq0$ and $m\geq1$. Then $\operatorname{supp} u(t) \subset \operatorname{supp} u_0$ a.e. in $(0,T]$.
\end{theorem}
{\bf Proof of Theorem \ref{nochange}.}
Let $\theta :\overline \Omega   \to \mathbb{R}$, which satisfies the conditions
\begin{align}
  \operatorname{supp} \;\theta  \subset \overline {\Omega \backslash \operatorname{supp} \;{u_0}} , \label{condi1} \\
  \theta  \in W_0^{1,\infty }\left( \Omega  \right).\label{condi2}
\end{align}
For example, we can take $\theta (x)= \min \left\{ {\dfrac{1}{\sigma }\operatorname{dist} \left( {x,\operatorname{supp}{u_0} \cup \partial \Omega } \right),1} \right\}$, $\sigma  \leqslant 1$. Observe that \eqref{condi1} implies $\theta  \cdot {u_0} \equiv 0$ on $\Omega$.
Taking $\phi : = \frac{\theta }{{u + \varepsilon }}$ as a test function in \eqref{weak form}, we obtain
\[\iint_{{\Omega _t}} {\frac{{\partial u}}{{\partial t}}}\frac{\theta }{{u + \varepsilon }}dxdt +  {\iint_{{\Omega _t}} {{{\left| {Du} \right|}^{p(x) - 2}}Du \cdot D\left( {{u^m}\frac{\theta }{{u + \varepsilon }}} \right)}dxdt}  = 0.\]
and we calculate
\begin{align*}
 & \int_\Omega  {\ln \left( {u(t) + \varepsilon } \right)\theta } dx - \int_\Omega  {\ln \left( {{u_0} + \varepsilon } \right)\theta } dx
  \\+ & {\iint_{{\Omega _t}} {\theta {{\left| {Du} \right|}^{p(x)}} \cdot \frac{{(m - 1){u^m} + \varepsilon m{u^{m - 1}}}}{{{{\left( {u + \varepsilon } \right)}^2}}}}dxdt}
  \\=&  {\iint_{{\Omega _t}} {\frac{{{u^m}}}{{u + \varepsilon }}{{\left| {Du} \right|}^{p(x) - 2}}Du \cdot D\theta }dxdt} .
\end{align*}
According to
\[\int_\Omega  {{\chi _{\left\{ {\operatorname{supp} \;\theta } \right\}}}\left( {\ln \left( {u(t) + \varepsilon } \right){-}\ln \left( \varepsilon  \right)} \right)\theta } dx \leqslant C,\]
for every $\delta$ sufficiently small, we have
\[\int_\Omega  {{\chi _{\left\{ {u(x,t) > \delta } \right\} \cap \{ \theta  = 1\} }}\left( {\ln \left( {u(t) + \varepsilon } \right){-}\ln \left( \varepsilon  \right)} \right)} dx \leqslant C,\]
where $C$ is independent of $\varepsilon$. We therefore conclude that
\begin{align*}
\text{measure}{\left\{ {(x,t) \in \{ \theta  = 1\}  \times \{ t\} \left| {u(x,t) > \delta } \right.} \right\}} = 0,
\text{ for a.e. } t\in (0,T),
\end{align*}
which implies the claim.$\hfill\blacksquare$
\begin{remark}
If $0<m<1$, Theorem \ref{nochange} no longer holds in general. For example, $p
(x)\equiv 2$, the problem \eqref{fu1}
has a Barenblatt solution in the following form:
\[{B_m }(x,t) = {(t + {t_0})^{ - \gamma}}{\left( {{{\left( {1 - \frac{m\gamma}{{2 N}} \frac{{\left| x \right|}^2}{{\left( {t + {t_0}} \right)}^{1-m\gamma}}} \right)}_ + }} \right)^{\frac{1}{m}}},\]
where $\gamma = \frac{N}{mN+2-2m}$, $N$ denotes the dimension of the spatial space. By calculation, one has
\begin{align*}
B_m^m\Delta {B_m} =  - \frac{\gamma }{{t + {t_0}}}{B_m}(x,t) + \frac{{\gamma (1 - m\gamma )}}{{2N}}\frac{{{{\left| x \right|}^2}}}{{{{(t + {t_0})}^2}}}B_m^{1 - m}(x,t) = \frac{{\partial {B_m}}}{{\partial t}}.
\end{align*}
We observe that $\operatorname{supp} \;{u_0} \subsetneqq \operatorname{supp} \;u(t)$ for $t>0$.
\end{remark}


\section*{Acknowledgement}
This work is partially supported by the National Natural Science
Foundation of China (11871133, 11971131, U1637208, 61873071, 51476047, 12171123), the Guangdong
Basic and Applied Basic Research Foundation (2020B1515310010, 2020B1515310006), the Fudament Research Funds for the Gentral Universities (Grant No.~NSRIF.~2020081), and the Natural Sciences Foundation of Heilongjiang Province (LH2020A004).

\bibliography{mybibfile}

\begin{thebibliography}{10}
\expandafter\ifx\csname url\endcsname\relax
  \def\url#1{\texttt{#1}}\fi
\expandafter\ifx\csname urlprefix\endcsname\relax\def\urlprefix{URL }\fi
\expandafter\ifx\csname href\endcsname\relax
  \def\href#1#2{#2} \def\path#1{#1}\fi

\bibitem{MR2811763}
S.~Antontsev, S.~Shmarev, Parabolic equations with double variable
  nonlinearities, Math. Comput. Simulation 81~(10) (2011) 2018--2032.

\bibitem{MR3013415}
S.~Antontsev, S.~Shmarev, Doubly degenerate parabolic equations with variable
  nonlinearity {I}: {E}xistence of bounded strong solutions, Adv. Differential
  Equations 17~(11-12) (2012) 1181--1212.

\bibitem{MR2438319}
M.~Y. Chen, Blow-up for doubly degenerate nonlinear parabolic equations, Acta
  Math. Sin. (Engl. Ser.) 24~(9) (2008) 1525--1532.

\bibitem{MR3130539}
S.~N. Antontsev, S.~Shmarev, Doubly degenerate parabolic equations with
  variable nonlinearity {II}: {B}low-up and extinction in a finite time,
  Nonlinear Anal. 95 (2014) 483--498.

\bibitem{MR965742}
M.~Bertsch, R.~Dal~Passo, M.~Ughi, Discontinuous ``viscosity'' solutions of a
  degenerate parabolic equation, Trans. Amer. Math. Soc. 320~(2) (1990)
  779--798.

\bibitem{MR1174811}
M.~Bertsch, R.~Dal~Passo, M.~Ughi, Nonuniqueness of solutions of a degenerate
  parabolic equation, Ann. Mat. Pura Appl. (4) 161 (1992) 57--81.

\bibitem{MR1044287}
M.~Bertsch, M.~Ughi, Positivity properties of viscosity solutions of a
  degenerate parabolic equation, Nonlinear Anal. 14~(7) (1990) 571--592.

\bibitem{MR853975}
A.~Friedman, B.~McLeod, Blow-up of solutions of nonlinear degenerate parabolic
  equations, Arch. Rational Mech. Anal. 96~(1) (1986) 55--80.

\bibitem{MR707172}
L.~J.~S. Allen, Persistence and extinction in single-species reaction-diffusion
  models, Bull. Math. Biol. 45~(2) (1983) 209--227.

\bibitem{MR859613}
M.~Ughi, A degenerate parabolic equation modelling the spread of an epidemic,
  Ann. Mat. Pura Appl. (4) 143 (1986) 385--400.

\bibitem{MR2902844}
J.~Yin, C.~Jin, Critical exponents of a doubly degenerate non-divergent
  parabolic equation with interior and boundary sources, Math. Nachr. 285~(5-6)
  (2012) 758--777.

\bibitem{MR2113162}
W.~Zhou, Z.~Wu, Some results on a class of degenerate parabolic equations not
  in divergence form, Nonlinear Anal. 60~(5) (2005) 863--886.

\bibitem{MR3635371}
C.~Jin, J.~Yin, Asymptotic behavior of solutions for a doubly degenerate
  parabolic non-divergence form equation, Rocky Mountain J. Math. 47~(2) (2017)
  479--510.

\bibitem{MR4107096}
J.~Shao, Z.~Guo, X.~Shan, C.~Zhang, B.~Wu, A new non-divergence diffusion
  equation with variable exponent for multiplicative noise removal, Nonlinear
  Anal. Real World Appl. 56 (2020) 103166, 15.

\bibitem{MR901093}
R.~Dal~Passo, S.~Luckhaus, A degenerate diffusion equation not in divergence
  form, in: Nonlinear parabolic equations: qualitative properties of solutions
  ({R}ome, 1985), Vol. 149 of Pitman Res. Notes Math. Ser., Longman Sci. Tech.,
  Harlow, 1987, pp. 72--76.

\bibitem{MR2171907}
W.~Zhou, Z.~Wu, Existence and nonuniqueness of weak solutions of the
  initial-boundary value problem for {$u_t=u^\sigma{\rm div}(|\nabla
  u|^{p-2}\nabla u)$}, Northeast. Math. J. 21~(2) (2005) 189--206.

\bibitem{MR1134951}
O.~Kov\'{a}\v{c}ik, J.~R\'{a}kosn\'{\i}k, On spaces {$L^{p(x)}$} and
  {$W^{k,p(x)}$}, Czechoslovak Math. J. 41(116)~(4) (1991) 592--618.

\bibitem{MR2790542}
L.~Diening, P.~Harjulehto, P.~H\"{a}st\"{o}, M.~Ruzicka, Lebesgue and {S}obolev
  spaces with variable exponents, Vol. 2017 of Lecture Notes in Mathematics,
  Springer, Heidelberg, 2011.

\bibitem{MR3328376}
S.~Antontsev, S.~Shmarev, Evolution {PDE}s with nonstandard growth conditions,
  Vol.~4 of Atlantis Studies in Differential Equations, Atlantis Press, Paris,
  2015, existence, uniqueness, localization, blow-up.

\end{thebibliography}

\clearpage
\appendix
\section{Proof of Lemmas}

{\bf{Proof of Lemma \ref{pmin}}}: If $2\leq p< \infty$, by the rearrangement inequality, we have
\begin{align}
& \left( {{{\left| \xi  \right|}^{p - 2}}\xi  - {{\left| \eta  \right|}^{p - 2}}\eta } \right) \cdot \left( {\xi  - \eta } \right) \nonumber
 \\  = &{\left| \xi  \right|^p} + {\left| \eta  \right|^p} - \left( {{{\left| \xi  \right|}^{p - 2}} + {{\left| \eta  \right|}^{p - 2}}} \right)\xi  \cdot \eta \nonumber
 \\  \geqslant& \left( {{{\left| \xi  \right|}^{p - 2}} + {{\left| \eta  \right|}^{p - 2}}} \right) \cdot \frac{{{{\left| \xi  \right|}^2} + {{\left| \eta  \right|}^2}}}{2} - \left(   {{{\left| \xi  \right|}^{p - 2}} + {{\left| \eta  \right|}^{p - 2}}} \right)\xi  \cdot \eta \nonumber
  \\ =& \frac{1}{2}\left( {{{\left| \xi  \right|}^{p - 2}} + {{\left| \eta  \right|}^{p - 2}}} \right){\left| {\xi  - \eta } \right|^2}. \label{le3.1}
\end{align}
For the case $2\leq p<3$, one has
\begin{align}
{\left| \xi  \right|^{p - 2}} + {\left| \eta  \right|^{p - 2}} \geqslant {\left( {\left| \xi  \right| + \left| \eta  \right|} \right)^{p - 2}} \geqslant {\left| {\xi  - \eta } \right|^{p - 2}}. \label{le3.2}
\end{align}
For the case $3\leq p$, by the convexity of $|\cdot|^{p-2}$, one has
\begin{align}
\frac{1}{2}{\left| \xi  \right|^{p - 2}}{\text{ + }}\frac{1}{2}{\left| \eta  \right|^{p - 2}} \geqslant {\left| {\frac{{\xi  - \eta }}{2}} \right|^{p - 2}}. \label{le3.3}
\end{align}
As a consequence of \eqref{le3.1}, \eqref{le3.2} and \eqref{le3.3}, we obtain $\rm{(i)}$.

If $\;1\leq p<2$. Assume that $\forall \;\theta \in [0,1]$, $\theta \xi +(1-\theta)\eta \ne 0$.
Using Cauchy's mean value theorem, we have
\begin{align}
  \left( {{{\left| \xi  \right|}^{p - 2}}\xi  - {{\left| \eta  \right|}^{p - 2}}\eta } \right) \cdot \left( {\xi  - \eta } \right)
=\int_0^1 {\left( {\xi  - \eta ,A(s)(\xi  - \eta )} \right)} \mathrm{d}s, \label{le3.4}
\end{align}
where
\begin{align*}
&A(s) = {\left( {{a_{ij}}(s)} \right)_{n \times n}},
\\&{a_{ij}}(s) = {\left| x \right|^{p - 2}}\left( {{\delta _{ij}} + (p - 2)\frac{{{x_i}{x_j}}}{{|x{|^2}}}} \right),
\\&x = ({x_i}) = \eta  + s(\xi  - \eta ) \in {\mathbb{R}^n},
\end{align*}
and $\delta_{ij}$ is Kronecker delta function.
\\On the other hand, for any $z\in\mathbb{R}^n$, it follows that
\begin{align}
  &\left( {z,A(s)z} \right) \nonumber\\
   = &{\left| x \right|^{p - 2}}\left( {|z{|^2} + (p - 2)\frac{{|x \cdot z{|^2}}}{{|x{|^2}}}} \right)  \nonumber\\
   \geqslant& {\left| x \right|^{p - 2}}\left( {|z{|^2} + (p - 2)\frac{{{{(|x||z|)}^2}}}{{|x{|^2}}}} \right)  \nonumber\\
   = &(p - 1){\left| x \right|^{p - 2}}|z{|^2}. \label{le3.5}
\end{align}
Then, by the convexity of $|\cdot|^p$, we get
\begin{align}
{\left| x \right|^p} \leqslant s{\left| \xi  \right|^p} + (1 - s){\left| \eta  \right|^p} \leqslant {\left| \xi  \right|^p} + {\left| \eta  \right|^p}. \label{le3.6}
\end{align}
Thus, by \eqref{le3.4}, \eqref{le3.5} and \eqref{le3.6}, we obtain
\begin{align}
\left( {{{\left| \xi  \right|}^{p - 2}}\xi  - {{\left| \eta  \right|}^{p - 2}}\eta } \right) \cdot \left( {\xi  - \eta } \right) \geq (p - 1){\left( {{{\left| \xi  \right|}^p} + {{\left| \eta  \right|}^p}} \right)^{\frac{{p - 2}}{p}}}{\left| {\xi  - \eta } \right|^2}. \label{le3.7}
\end{align}
If $\;1\leq p<2$, $\exists \;\theta \in [0,1]$ such that $\theta \xi +(1-\theta)\eta = 0$, we shall only prove
\begin{align}
({k^{p - 1}} - 1)(k - 1) \geqslant (p - 1){(k - 1)^2}{({k^p} + 1)^{\frac{{p - 2}}{p}}},\;\forall k \geqslant 0. \label{le3.8}
\end{align}
The inequality \eqref{le3.8} is based on similar arguments with \eqref{le3.7}. Collecting all these facts, we complete the proof of $\rm{(ii)}$. $\hfill\blacksquare$

{\bf{Proof of Lemma \ref{holder} }}:
Let us denote ${\left\| f \right\|_{q( \cdot ),\Omega }}= \lambda$, ${\left\| g \right\|_{q'( \cdot ),\Omega }}=\mu$ and assume that $\lambda\mu \ne 0$. By Young's inequality, one has for a.e. $x\in \Omega$,
\begin{align}
&\left| {f(x)g(x)} \right| \nonumber
\\=& \lambda \mu \left| {\frac{{f(x)}}{\lambda }} \right|\left| {\frac{{g(x)}}{\mu }} \right| \nonumber
\\\leqslant &\lambda \mu \left( {\frac{1}{{q(x)}}{{\left| {\frac{{f(x)}}{\lambda }} \right|}^{q(x)}} + \frac{1}{{q'(x)}}{{\left| {\frac{{g(x)}}{\mu }} \right|}^{q'(x)}}} \right) \nonumber
\\  \leqslant & \lambda \mu \left( {{{\left| {\frac{{f(x)}}{\lambda }} \right|}^{q(x)}} + {{\left| {\frac{{g(x)}}{\mu }} \right|}^{q'(x)}}} \right). \label{hi1}
\end{align}
On the other hand, by the definition of Luxemburg norm and monotone convergence theorem, one has
\begin{align}
{\varrho _{q( \cdot )}}(\frac{f}{\lambda }) \leqslant 1,\quad {\varrho _{q'( \cdot )}}(\frac{g}{\mu }) \leqslant 1. \label{hi2}
\end{align}
Integrating \eqref{hi1} over $\Omega$ and applying \eqref{hi2}, we have
\begin{align*}
&\int_\Omega  {\left| {f(x)g(x)} \right|} \mathrm{d}x
\\ \leqslant &\lambda \mu \left( {{\varrho _{q( \cdot )}}(\frac{f}{\lambda }) + {\varrho _{q'( \cdot )}}(\frac{g}{\mu })} \right)
\\   \leqslant & 2{\left\| f \right\|_{q( \cdot ),\Omega }}{\left\| g \right\|_{q'( \cdot ),\Omega }}.
\end{align*}
In the case $\lambda\mu = 0$, the inequality is trivial.   $\hfill\blacksquare$

{\bf{Proof of Lemma \ref{lest}}}:
According to $\rm{(ii)}$ of Lemma \ref{pmin}, we have
\begin{align}
  &\int_\Omega  {\left( {{{\left| u \right|}^{p(x) - 2}}u - {{\left| v \right|}^{p(x) - 2}}v} \right) \cdot \left( {u - v} \right)\mathrm{d}x} \nonumber\\
   \geqslant &({p^ - } - 1)\int_\Omega  {{{\left| {u - v} \right|}^2} \cdot {{(|u{|^{p(x)}} + |v{|^{p(x)}})}^{\frac{{p(x) - 2}}{{p(x)}}}}\mathrm{d}x}. \label{le4.1}
\end{align}
On the other hands, according to Lemma \ref{holder}, we obtain
 \begin{align}
&\int_\Omega  {{{\left| {u - v} \right|}^{p(x)}}\mathrm{d}x}  \cdot {\left( 2{{{\left\| {{{(|u{|^{p( \cdot )}} + |v{|^{p( \cdot )}})}^{\frac{{2 - p( \cdot )}}{2}}}} \right\|}_{\frac{2}{{2 - p( \cdot )}},\Omega }}} \right)^{ - 1}} \nonumber
\\ \leqslant & {\left\| {|u - v{|^{p( \cdot )}}{{(|u{|^{p( \cdot )}} + |v{|^{p( \cdot )}})}^{\frac{{p( \cdot ) - 2}}{2}}}} \right\|_{\frac{2}{{p( \cdot )}},\Omega }}. \label{le4.2}
 \end{align}
 Combining \eqref{le4.1}, \eqref{le4.2} and Lemma \ref{vels}, the conclusion follows.$\hfill\blacksquare$

\section{Weak Solution of Regularized Problem}
In order to obtain the weak solution of regularized problem, we discuss the conditions in three cases.
\begin{itemize}
\item
For the case $0<m<1$, let
\begin{gather*}
  v = \Phi (u) = \frac{{{u^{1 - m}}}}{{1 - m}},\\
  u = \Psi (v) = {\left( {(1 - m)v} \right)^{\frac{1}{{1 - m}}}}.
\end{gather*}
\item
For the case $m> 1$, let
\begin{gather*}
  v = \Phi (u) = \frac{{{u^{1 - m}}}}{{m-1}},\\
  u = \Psi (v) = {\left( {(m-1)v} \right)^{\frac{1}{{1 - m}}}}.
\end{gather*}
\item
For the case $m=1$, let
\begin{gather*}
  v = \Phi (u) = \ln u,\\
  u = \Psi (v) = {e^v}.
\end{gather*}
\end{itemize}
Then the problem \eqref{reg1} is translated into parabolic equations in divergence form as follows:
\begin{numcases}{}
  {v_t} =  {{\rm div}\left( {{{\left| {\Psi '(v)} \right|}^{{p}(x) - 1}}{{\left| {{D}v} \right|}^{{p}(x) - 2}}{D}v} \right)}  \quad &\text{in } $\Omega_T$, \nonumber\\
  v(x,t) = \Phi (\varepsilon ) \quad &\text{on } $\Gamma$,\label{reg5} \\
  v(x,0) = \Phi ({u_0} + \varepsilon )\quad &\text{in } $\Omega$. \nonumber
\end{numcases}

\begin{definition}
A function $v(x,t)$ is called a weak solution of parabolic problem \eqref{reg5} provided that
\begin{itemize}
\item
$v \in \mathcal{U}\left( {{\Omega _T}} \right) \cap {L^\infty }\left( {{\Omega _T}} \right),\;\;v_t \in \mathcal{U}'\left( {{\Omega _T}} \right) .$
\item
For every $\varphi \in C_0^1 (\Omega_T)$,
\begin{align}
\iint_{{\Omega _T}} {{v_t}}\varphi \mathrm{d}x\mathrm{d}t + \iint_{{\Omega _T}} {{{{\left| {\Psi '(v)} \right|}^{{p}(x) - 1}}{{\left| {{D}v} \right|}^{p(x) - 2}}{D}v \cdot {D}\varphi } }\mathrm{d}x\mathrm{d}t = 0.
\label{reg1 weak form}
\end{align}
\item
The following equations hold in the sense of trace:
\begin{align*}
v(x,t) = \Phi (\varepsilon ) \qquad &\text{on } \Gamma,
\\
v(x,0) = \Phi ({u_0} + \varepsilon ) \qquad &\text{in } \Omega.
\end{align*}
\end{itemize}
\end{definition}

Denote
\[K = \mathrm{ess} \mathop { \sup }\limits_{x \in \Omega } \;{u_0}(x),\]
\[{A_{\varepsilon ,K}} = \max \left\{ {{\varepsilon ^m},\min \left\{ {{ \left| {\Psi '(v)} \right|  },{{(K + \varepsilon )}^m}} \right\}} \right\}.\]
Then we consider the regular equations as follows:
\begin{numcases}{}
 {v_t} = {\rm div} {\left( {A_{\varepsilon ,K}^{{p}(x) - 1}{{\left| {{D}v} \right|}^{{p}(x) - 2}}{D}v} \right)}, \quad &\text{in } $\Omega_T$, \nonumber \\
  v(x,t) = \Phi (\varepsilon ) \quad &\text{on } $\Gamma$,\label{reg7} \\
  v(x,0) = \Phi ({u_0} + \varepsilon )\quad &\text{in } $\Omega$. \nonumber
\end{numcases}

\begin{definition}
A function $v(x,t)$ is called a weak solution of regularized problem \eqref{reg7} provided that
\begin{itemize}
\item
$v \in \mathcal{U}\left( {{\Omega _T}} \right) \cap {L^\infty }\left( {{\Omega _T}} \right),\;\;v_t \in \mathcal{U}'\left( {{\Omega _T}} \right) .$
\item
For every $\varphi \in C_0^1 (\Omega_T)$,
\begin{align}
\iint_{{\Omega _T}} {{v_t}}\varphi \mathrm{d}x\mathrm{d}t + \iint_{{\Omega _T}} { {A_{\varepsilon ,K}^{{p}(x) - 1}{{\left| {{D}v} \right|}^{p(x) - 2}}{D}v \cdot {D}\varphi } }\mathrm{d}x\mathrm{d}t = 0.
\label{ reg2 weak form}
\end{align}
\item
The following equations hold in the sense of trace:
\begin{align*}
v(x,t) = \Phi (\varepsilon ) \qquad &\text{on } \Gamma,
\\
v(x,0) = \Phi ({u_0} + \varepsilon ) \qquad &\text{in } \Omega.
\end{align*}
\end{itemize}
\end{definition}

\begin{remark}{\label{vir}}
We consider the equation of $w(x,t)=v(x,t)-\Phi (\varepsilon)$, then $w(x,t) = 0,\;(x, t) \in \Gamma$. By virtue of \cite[Theorem 4.1]{MR3328376}, there exists a weak solution $w(x,t)$. So the regular problem \eqref{reg7} admits a weak solution.
\end{remark}

\begin{proposition} \label{max}
For the case $0<m\leq 1$, the weak solution of the regular problem \eqref{reg7} satisfies
\[\Phi \left( \varepsilon  \right) \leqslant v(x,t) \leqslant \Phi \left( {K + \varepsilon } \right).\]
For the case $m>1$, the weak solution of the regular problem \eqref{reg7} satisfies
\[0<\Phi \left( {K + \varepsilon } \right) \leqslant v(x,t) \leqslant \Phi \left( \varepsilon  \right).\]
\end{proposition}
{\bf{Proof }}:
For the case $0<m\leq 1$, multiplying the equation \eqref{reg7} by ${\left( {v - M} \right)_ + }$, and integrating over $\Omega_s$, we have
\[{\frac{1}{2}} \iint_{\Omega_s} \frac{\partial }{{\partial t}}\left( {v - M} \right)_ + ^2\mathrm{d}x\mathrm{d}t = - \iint_{{\Omega _s}} { {A_{\varepsilon ,K}^{{p}(x) - 1}{{\left| {{D}{{\left( {v - M} \right)}_ + }} \right|}^{{p}(x)}}} }\mathrm{d}x\mathrm{d}t \leqslant 0,\]
where $M>0$ is a constant which will be determined later.

Therefore,
\[\int_\Omega  {\left( {v(x,s) - M} \right)_ + ^2} \mathrm{d}x \leqslant \int_\Omega  {\left( {v(x,0) - M} \right)_ + ^2} \mathrm{d}x.\]
Due to $v(x,0)\leq \Phi (K + \varepsilon )$ and the arbitrariness of $s$, choosing $M=\Phi \left( {K + \varepsilon } \right)$, we have $v(x,t) \leqslant \Phi \left( {K + \varepsilon } \right)$ a.e. in $\Omega_T$. Similarly multiplying the equation \eqref{reg7} by ${\left( {v - N} \right)_ - }$, choosing $N=\Phi \left( \varepsilon  \right)$, we have $\Phi \left( \varepsilon  \right) \leqslant v(x,t)$.

In a similar way, we can get the conclusion for the case $1<m< 2$.$\hfill\blacksquare$

Based on Proposition \ref{max}, we know that $0 < {\varepsilon ^m} \leqslant \left| {\Psi '(v)} \right| \leqslant {\left( {K + \varepsilon } \right)^m}$.
Thus, the weak solution of \eqref{reg7} is the weak solution of \eqref{reg5}. The following Corollary follows.

\begin{corollary}
Assume that $m>0$ and $p(x)$ is log-H\"older continuous which satisfies \eqref{px}, $u_0 \in L^{\infty}(\Omega)$ and $u_0 \geq0$, the problem \eqref{reg5} admits a weak solution.
\end{corollary}

{\bf{Proof of Proposition \ref{reex}}}:
Since that $p_i(x)$ is log-H\"older continuous, and $C_0^1 (\Omega_T)$ is dense in ${\mathcal{U}_0}\left( {{\Omega _T}} \right)$, \eqref{reg1 weak form} holds true also for all
$\varphi \in {\mathcal{U}_0}\left( {{\Omega _T}} \right)$.
\\For any $\phi \in C_0^1 (\Omega_T)$, taking $\varphi  = \left| {\Psi '(v)} \right|\phi $ in \eqref{reg1 weak form},
we have
\[\iint_{{\Omega _T}} {\Psi {{(v)}_t}}\phi \mathrm{d}x\mathrm{d}t + \iint_{{\Omega _T}} { {{{\left| {{D}\Psi (v)} \right|}^{p(x) - 2}}{D}\Psi (v) \cdot {D}\left( {\left| {\Psi '(v)} \right|\phi } \right)} }\mathrm{d}x\mathrm{d}t = 0.\]
In fact of ${u^m} = \left| {\Psi '(v)} \right|$, \eqref{reg weak form} holds for $u={\Psi (v)}$.
\\Since
\begin{align*}
v(x,t) = \Phi (\varepsilon ) \qquad &\text{on } \Gamma,
\\
v(x,0) = \Phi ({u_0} + \varepsilon ) \qquad &\text{in } \Omega ,
\end{align*}
we have \eqref{trace3}, \eqref{trace4} in the sense of trace, then $u={\Psi (v)}$ is a weak solution of \eqref{reg1}.  $\hfill\blacksquare$

{\bf{Proof of Proposition \ref{bj}}}: Assume $u_1$ and $u_2$ are the solutions of the equation which correspond to $\varepsilon_1$ and $\varepsilon_2$ respectively, and $\varepsilon_1 \leq \varepsilon_2$. Choosing $u_1^{ - m}H\left( {{u_1} - {u_2}} \right)$ and
$u_2^{ - m}H\left( {{u_1} - {u_2}} \right)$ as the test function, where
\[H(t) = \left\{ \begin{gathered}
  1,\;\;t > 0,  \\
  0,\;\;t \leqslant 0,
\end{gathered}  \right.\]
then we have
\begin{align*}
&\iint_{{\Omega_T}} \left ( u_1^{ - m}{u_{1t}} - u_2^{ - m}{u_{2t}} \right ) H\left( {{u_1} - {u_2}} \right)\mathrm{d}x\mathrm{d}t  \\+
&\iint_{{\Omega _T}} {\delta \left( {{u_1} - {u_2}} \right)\cdot\sum\limits_{i = 1}^n {\left( {{{\left| {{D}{u_1}} \right|}^{p(x) - 2}}{D_i}{u_1} - {{\left| {{D}{u_2}} \right|}^{p(x) - 2}}{D}{u_2}} \right) \cdot {D}\left( {u_1-u_2}\right)} }\mathrm{d}x\mathrm{d}t = 0.
\end{align*}

Now, for the case $m\neq1$,
\begin{align*}
&\frac{1}{{1 - m}}\int_{\left\{ {x \in \Omega ;{u_1} \geqslant {u_2}} \right\}} {u_1^{1 - m}\left( {x,T} \right) - u_2^{1 - m}\left( {x,T} \right)} \mathrm{d}x \\
\leqslant & \frac{1}{{1 - m}}\int_{\left\{ {x \in \Omega ;{u_{10}} \geqslant {u_{20}}} \right\}} {u_1^{1 - m}\left( {x,0} \right) - u_2^{1 - m}\left( {x,0} \right)} \mathrm{d}x = 0.
\end{align*}
Likewise, for the case $m=1$,
\begin{align*}
&\int_{\left\{ {x \in \Omega ;{u_1} \geqslant {u_2}} \right\}} {\ln ({u_1}\left( {x,T}) \right) - \ln({u_2}\left( {x,T}) \right)} \mathrm{d}x
\\ \leqslant &\int_{\left\{ {x \in \Omega ;{u_{10}} \geqslant {u_{20}}} \right\}} {\ln({u_1}\left( {x,0}) \right) - \ln ({u_2}\left( {x,0} )\right)} \mathrm{d}x = 0.
\end{align*}
Therefore, $u_1 \leq u_2$ a.e. in $\Omega_T$.$\hfill\blacksquare$
\begin{remark}
In fact, we can complete the proof through a process of approximation; that is, we can choose $H_\epsilon (t)$ instead of $H(t)$, where
\[{H_\epsilon}(t) = \int_0^t {{h_\epsilon }(s)\mathrm{d}s} ,\;\;{h_\epsilon }(t) = \frac{2}{\epsilon }{\left( {1 - \frac{{\left| s \right|}}{\epsilon }} \right)_ + },\]
and then let $\epsilon\rightarrow 0$.
\end{remark}


\section{Proof of Proposition \ref{Estimates}}
In order to obtain the weak solution of problem \eqref{fu1}, some apriori estimates are also necessary.

Assume that $0\leq u_0 \in \mathcal {V}(\Omega ) \cap L^{\infty}(\Omega)$ and $p_i(x)$ are log-H\"older continuous functions which satisfy \eqref{px}.

Due to $u = \Psi (v)= \Phi^{-1}(v)$, choosing $\varepsilon$ small enough, we know from Proposition \ref{max} that
\begin{align}
0< \varepsilon \leq u_\varepsilon \leq K+\varepsilon \leq K+1.  \label{est0}
\end{align}

Multiplying the equation \eqref{reg1} by $u_\varepsilon ^{ - m}\dfrac{{\partial {u_\varepsilon }}}{{\partial t}}$, and integrating over $\Omega_T$, we have
\[\iint_{{\Omega _T}} {u_\varepsilon ^{ - m}}{\left( {\frac{{\partial {u_\varepsilon }}}{{\partial t}}} \right)^2}\mathrm{d}x\mathrm{d}t +  {\int_\Omega  {\frac{1}{{{p}(x)}}{{\left| {{D}{u_\varepsilon }} \right|}^{p(x)}}} \mathrm{d}x}  = {\int_\Omega  {\frac{1}{{{p}(x)}}{{\left| {{D}{u_0}} \right|}^{p(x)}}} \mathrm{d}x} ,\]
then
\begin{align}
  \iint_{{\Omega_T}} {u_\varepsilon ^{ - m}}{\left( {\frac{{\partial {u_\varepsilon }}}{{\partial t}}} \right)^2}\mathrm{d}x\mathrm{d}t \leqslant C, \label{est1} \\
   {\int_\Omega  {\frac{1}{{{p}(x)}}{{\left| {{D}{u_\varepsilon }} \right|}^{p(x)}}} \mathrm{d}x}  \leqslant C,\label{est2}
\end{align}
where $C$ is a constant independent of $\varepsilon$.

According to \eqref{est0} and \eqref{est1}, we have
\begin{align}
\iint_{{\Omega_T}} {{{\left( {\frac{{\partial {u_\varepsilon }}}{{\partial t}}} \right)}^\alpha }}\mathrm{d}x\mathrm{d}t
\leq  (K+1)^m \left( {\iint_{{\Omega_T}} {u_\varepsilon ^{ - m}}{{\left( {\frac{{\partial {u_\varepsilon }}}{{\partial t}}} \right)}^2}\mathrm{d}x\mathrm{d}t} \right) \leq C, \label{est3}
\end{align}
where $C$ is a constant independent of $\varepsilon$.

\begin{remark}
Actually, the process above can be completed by apriori estimates of the regular equations \eqref{reg7}. Denote the solution of the regular problem \eqref{reg7} as $v^\varepsilon$. Multiplying the equation \eqref{reg7} by $\dfrac{{\partial \Psi (v^\varepsilon)}}{{\partial t}}$, and integrating over $\Omega_T$, we have
\[\iint_{{\Omega _T}} {\left| {\Psi '({v^\varepsilon })} \right|}{\left( {v_t^\varepsilon } \right)^2}\mathrm{d}x\mathrm{d}t +  {\iint_{{\Omega _T}} {{{\left| {{D}\Psi ({v^\varepsilon })} \right|}^{{p}(x) - 2}}{D}\Psi ({v^\varepsilon }) \cdot \frac{\partial }{{\partial t}}\left( {{D}\Psi ({v^\varepsilon })} \right)}\mathrm{d}x\mathrm{d}t}  = 0.\]
Therefore,
\[\iint_{{\Omega _T}} {\left| {\Psi '({v^\varepsilon })} \right|}{\left( {v_t^\varepsilon } \right)^2}\mathrm{d}x\mathrm{d}t +  {\int_\Omega  {\frac{1}{{{p}(x)}}{{\left| {{D}\Psi ({v^\varepsilon })} \right|}^{{p}(x)}}} \mathrm{d}x}  =  {\int_\Omega  {\frac{1}{{{p}(x)}}{{\left| {{D}{u_0}} \right|}^{{p}(x)}}} \mathrm{d}x} .\]
Then \eqref{est1}--\eqref{est3} follows.
\end{remark}

\end{document}